\documentclass[12pt,righttagstyle]{amsart}

\usepackage{amssymb,latexsym}

\newcommand{\epf}{ $\Box$\medskip}

\newtheorem{thm}{Theorem}[section]

\newtheorem{lem}[thm]{Lemma}
\newtheorem{defn}[thm]{Definition}
\newtheorem{remark}[thm]{Remark}
\begin{document}

\title[norms inequalities for maximal operator]{Weighted norms inequalities for a maximal operator in some subspace of amalgams.}

\author[J. Feuto]{Justin Feuto}
\address{UFR DE MATHEMATIQUES ET INFORMATIQUE,
 Universit\'e de Cocody, 22BP 1194 Abidjan 22, R\'epublique de C\^ote d'Ivoire}
\email{{\tt justfeuto@yahoo.fr}}
\author[I. Fofana]{Ibrahim Fofana}
\address{UFR DE MATHEMATIQUES ET INFORMATIQUE,
 Universit\'e de Cocody, 22BP582 Abidjan 22, R\'epublique de C\^ote d'Ivoire}
\email{{\tt fofana\_ib\_math\_ab@yahoo.fr}}
\author[K. Koua]{Konin Koua}
\address{UFR DE MATHEMATIQUES ET INFORMATIQUE,
 Universit\'e de Cocody, 22BP582 Abidjan 22, R\'epublique de C\^ote d'Ivoire}
\email{{\tt kroubla@yahoo.fr }}

\thanks{}

\begin{abstract}
We give weighted norm inequalities for the maximal fractional operator $ \mathcal M_{q,\beta }$ of Hardy-Littlewood and the fractional integral $I_{\gamma}$. These inequalities are established between $\left( L^{q},L^{p}\right) ^{\alpha }(X,d,\mu )$ spaces (which are super spaces of Lebesgue spaces $L^{\alpha}(X,d,\mu)$, and subspaces of amalgams $(L^{q},L^{p})(X,d,\mu)$) and in the setting of space of homogeneous type $(X,d,\mu)$. The conditions on the weights are stated in terms of Orlicz norm.\\

\medskip 

{R\scshape\'esum\'e.} Nous donnons des in\'egalit\'es \`a poids pour l'op\'erateur maximal fractionnaire $\mathcal M_{q,\beta }$ de Hardy-Littlewood  et de
l'int\'{e}grale fractionnaire $I_{\gamma}$. Ces in\'egalit\'es sont \'etablies entre des espaces $\left( L^{q},L^{p}\right) ^{\alpha }(X,d,\mu )$ (qui sont des sur-espaces des espaces de Lebesgue $L^{\alpha}(X,d,\mu)$ et des sous-espaces des espaces d'amalgames $(L^{q},L^{p})(X,d,\mu)$) et dans le contexte des espaces de type homog\`{e}nes. Les conditions impos\'ees aux poids sont exprim\'ees en terme de norme d'Orlicz.

\end{abstract}

\maketitle

\section{Introduction}
Consider the fractional maximal operator $\mathfrak m_{q,\beta}$ ($1\leq q\leq\beta\leq\infty$) defined on $\mathbb R^{n}$ by 
\begin{equation}
\mathfrak m_{q,\beta}f(x)=\sup_{Q\in\mathcal Q:x\in Q}\left|Q\right|^{\frac{1}{\beta}-\frac{1}{q}}\left\|f\chi_{Q}\right\|_{q}
\end{equation}
where $\mathcal Q$ is the set of all cubes $Q$ of $\mathbb R^{n}$ with edges parallel to the coordinate axes, $\left|E\right|$ stands for the Lebesgue measure of the subset $E$ of $\mathbb R^{n}$ and $\left\|\cdot\right\|_{q}$ denotes the usual norm on the Lebesgue space $L^{q}(\mathbb R^{n},dx)$.
Weighted norm inequalities for $\mathfrak m_{1,\beta}$ have been extensively studied in the setting of Lebesgue, weak-Lebesgue and Morrey spaces (see \cite{3},\cite{8},\cite{N} and the references therein). The following result is contained in \cite{3}.
\begin{thm}\label{theoA} Assume that $1\leq q<\beta\leq\infty,$ $\frac{1}{t}=\frac{1}{q}-\frac{1}{\beta}$ and $v$ is a weight function satisfying
\begin{equation}
\sup_{Q\in\mathcal Q}\left|Q\right|^{\frac{1}{\beta}-1}\left\|v\chi_{_{Q}}\right\|_{t}\left\|v^{-1}\chi_{_{Q}}\right\|_{q'}<\infty.\ \ \ \ (\frac{1}{q'}+\frac{1}{q}=1)\label{AP}
\end{equation}
Then there exists a constant $C$ such that for any measurable function $f$
\begin{equation}
\left(\int_{\left\{x\in\mathbb R^{n}:\mathfrak m_{1,\beta}f(x)>\lambda\right\}}v(y)^{t}dy\right)^{\frac{1}{t}}\leq C\lambda^{-1}\left\|fv\right\|_{q}\ \ \ \lambda>0.
\end{equation}
\end{thm}

The spaces $(L^{q},\ell^{p})^{\alpha}(\mathbb R^{n})$ $(1\leq q\leq\alpha\leq p\leq\infty)$ have been defined in \cite{F2} as follows:\\
$\bullet\ I_{k}^{r}=\overset{n}{\underset{j=1}{\Pi }}\left[ k_{j}r,\left(
k_{j}+1\right) r\right[, \ \ \ k=(k_{i})_{1\leq i\leq n}\in\mathbb Z^{n},\ \ r>0$\\
$\bullet \ J_{x}^{r}=\overset{n}{\underset{j=1}{\Pi }}%
\left] x_{j}-\frac{r}{2}, x_{j}+\frac{r}{2}\right[,\ \ \ \ x=(x_{i})_{1\leq i\leq n}\in\mathbb R^{n},\ \ \ r>0$\\
$\bullet$ a Lebesgue measurable function $f$ belongs  to $(L^{q},\ell^{p})^{\alpha}(\mathbb R^{n})$ if $\left\| f\right\| _{q,p,\alpha }<\infty$, where
\begin{equation}
\left\| f\right\| _{q,p,\alpha }=\sup_{r>0}r^{n(\frac{1}{\alpha}-\frac{1}{q}}\ _{r}\left\|f\right\|_{q,p},
\end{equation}
and
\begin{equation}
 _{r}\left\|f\right\|_{q,p}=\left\{ 
\begin{array}{lll}
\left[\underset{k\in \mathbb{Z}^{n}}{\sum }\left( \left\| f\chi
_{_{I_{k}^{r}}}\right\| _{q}\right) ^{p}\right] ^{\frac{1}{p}} & \text{if }
& p< \infty \\ 
\underset{x\in \mathbb{R}^{n}}{\sup }\left\| f\chi _{_{J_{x}^{r}}}\right\|
_{q} & \text{if} & p=\infty%
\end{array}
\right..
\end{equation}
The $(L^{q},\ell^{p})^{\alpha}(\mathbb R^{n})$ have been introduced in connection with Fourier multiplier problems. But they are also linked to $L^{q}-L^{p}$ multiplier problems. We refer the readers to \cite{FKK} where space of Radon measures containing $(L^{1},\ell^{p})^{\alpha}(\mathbb R^{n})$ are considered. 

Notice that these spaces are subspaces of amalgams spaces introduced by Wiener and study by may authors (see the survey paper \cite{FS} of Fournier and Stewart and the references therein).

It has been proved in \cite{F1} that given $1\leq q\leq\alpha<\infty$, $\left\{(L^{q},\ell^{p})^{\alpha}(\mathbb R^{n})\right\}_{p\geq\alpha}$ is a monotone increasing family of Banach spaces, $(L^{q},\ell^{\alpha})^{\alpha}(\mathbb R^{n})=L^{\alpha}(\mathbb R^{n})$ and $(L^{q},\ell^{\infty})^{\alpha}(\mathbb R^{n})$ is clearly the classical Morrey space denoted by $L^{q,n(1-\frac{q}{\alpha})}(\mathbb R^{n})$ in \cite{8}. Moreover if $q<\alpha<p$ then the weak-$L^{\alpha}(\mathbb R^{n})$ space is embeded in $(L^{q},\ell^{p})^{\alpha}(\mathbb R^{n})$. Due to this remarkable link between the spaces $(L^{q},\ell^{p})^{\alpha}(\mathbb R^{n})$ and the Lebesgues ones, it is tempting to look for an extension of Theorem \ref{theoA} to the setting of $(L^{q},\ell^{p})^{\alpha}(\mathbb R^{n})$ space. The following result is contained in \cite{F3}.
\begin{thm}\label{theoB}
Assume that:\\
$\bullet \ 1\leq q\leq\alpha\leq p$ and  $0<\frac{1}{s}=\frac{1}{\alpha}-\frac{1}{\beta}$ \\
$\bullet\ q\leq q_{1}\leq\alpha_{1}\leq p_{1}$ and $0<\frac{1}{t}=\frac{1}{q_{1}}-\frac{1}{\beta}\leq\frac{1}{p_{1}}$,\\
$\bullet\ v$ is a weight function satisfying 
\begin{equation}
\sup_{Q\in\mathcal Q}\left|Q\right|^{\frac{1}{\beta}-\frac{1}{q}}\left\|v\chi_{Q}\right\|_{t}\left\|v^{-1}\chi_{Q}\right\|_{1/(\frac{1}{q}-\frac{1}{q_{1}})}<\infty.
\end{equation}
 Then there exists a real constant $C>0$ such that
\begin{equation}
\left(\int_{\left\{x\in\mathbb R^{n}:\mathfrak m_{1,\beta}f(x)>\lambda\right\}}v(y)^{t}dy\right)^{\frac{1}{t}}\leq C\lambda^{-1}\left\|fv\right\|_{q_{1},p_{1},\alpha_{1}}\left(\lambda^{-1}\left\|f\right\|_{q,\infty,\alpha}\right)^{s\left(\frac{1}{q_{1}}-\frac{1}{\alpha_{1}}\right)}
\end{equation}
for any real $\lambda>0$ and Lebesgue-measurable function $f$ on $\mathbb R^{n}$.
\end{thm}
It turns out that the $(L^{q},\ell^{p})^{\alpha}(\mathbb R^{n})$ setting is particularly adapted for the search of controls on Lebesgue norm of fractional maximal functions $\mathfrak m_{q,\beta}f$. Actually we have the following result which first part is a consequence of Theorem \ref{theoB}.

\begin{thm}(see \cite{F3})\label{theoC}
Assume that $1\leq q\leq\alpha\leq \beta$, and $\frac{1}{s}=\frac{1}{\alpha}-\frac{1}{\beta}$.
\begin{enumerate}
\item[(a)] If $\alpha \leq p$ and $\frac{1}{q}-\frac{1}{\beta}\leq\frac{1}{p}$ then  there is a real constant $C$ such that for all Lebesgue measurable function $f$ on $\mathbb R^{n}$,
\begin{equation}
\left\|\mathfrak m_{q,\beta}f\right\|^{\ast}_{s,\infty}\equiv\sup_{\lambda>0}\lambda\left|\left\{x\in\mathbb R^{n}:\mathfrak m_{q,\beta}f(x)>\lambda\right\}\right|^{\frac{1}{s}}\leq C\left\|f\right\|_{q,p,\alpha}.\label{inC1}
\end{equation}
\item [(b)] If $1\leq u\leq s\leq v$ then there is a real constant $C$ such that for any Lebesgue measurable function $f$ on $\mathbb R^{n},$ \begin{equation}
\left\|f\right\|_{q,p,\alpha}\leq C\left\|\mathfrak m_{q,\beta}f\right\|_{u,v,s}.\label{inC2}
\end{equation}
\end{enumerate}
\end{thm}
From inequality (\ref{inC2}) and the imbedding of the weak-$L^{s}(\mathbb R^{n})$ space into $(L^{u},\ell^{v})^{s}(\mathbb R^{n})$ for $u<s<v$, it follows that $f$ has its fractional maximal function $\mathfrak m_{q,\beta}f$ in a weak Lebesgue space only if itself belongs to some $(L^{q},\ell^{p})^{\alpha}(\mathbb R^{n})$.

Let  $X=(X,d,\mu)$ be a space of homogeneous type which is separable and satisfies a reversed doubling condition (see (\ref{revd}) in section 2 for definition). 
 
  For $1\leq q\leq\beta\leq\infty$ we set, for any $\mu$-measurable function $f$ on $X$,
  \begin{equation}
  \mathcal M_{q,\beta}f(x)=\sup_{B}\mu\left(B\right)^{\frac{1}{\beta}-\frac{1}{q}}\left\|f\chi_{B}\right\|_{q}\  \ \ x\in X,
  \end{equation}
  where the supremum is taken over all balls $B$ in $X$ containing $x$ and $\left\|\cdot\right\|_{q}$ denotes the norm on the Lebesgue space $L^{q}=L^{q}(X,d,\mu)$. As we can see, $\mathcal M_{q,\beta}$ is clearly a generalization of $\mathfrak m_{q,\beta}$. In the last decades, much work has been dedicated to obtain Morrey and Lebesgue norm inequalities for $\mathcal M_{q,\beta}$ and other operator of fractional maximal type on spaces of homogeneous type. We refer the reader to \cite{BS},\cite{6},\cite{EKM},\cite{1},\cite{2},\cite{sw},\cite{4} and the references therein.

As in the euclidean case, Lebesgue and Morrey spaces on homogeneous type spaces may be viewed as the end of a chain of Banach function spaces $(L^{q},L^{p})^{\alpha}(X)$ defined as follows: a $\mu$-measurable function $f$ represents an element of $\left( L^{q},L^{p}\right) ^{\alpha }(X)$ if  
\begin{equation}
\left\| f\right\| _{q,p,\alpha }=\sup_{r>0}\ _{r}\left\| f\right\| _{q,p,\alpha}<\infty 
\end{equation}
where
\begin{equation}
 _{r}\left\| f\right\| _{q,p,\alpha}=\left\{ 
\begin{array}{lll}
\left[ \int_{X}\left(\mu(B_{(y,r)})^{\frac{1}{\alpha}-\frac{1}{p}-\frac{1}{q}} \left\| f\chi _{_{B_{\left( y,r\right)
}}}\right\| _{q}\right) ^{p}d\mu (y)\right] ^{\frac{1}{p}} & \text{if} & 
p<\infty \\ 
\underset{y\in X}{\sup\text{ess}}\mu(B_{(y,r)})^{\frac{1}{\alpha}-\frac{1}{q}}\left\| f\chi _{_{B_{\left( y,r\right)
}}}\right\| _{q} & \text{if} & p=\infty%
\end{array}
\right. ,\label{namalgamX}
\end{equation}
The $\left( L^{q},L^{p}\right) ^{\alpha }(X)$ are generalizations of the $(L^{q},\ell^{p})^{\alpha}(\mathbb R^{n})$ which main properties  extend to them (see \cite{f-f}). 

In this paper we are interested in continuity properties of $\mathcal M_{q,\beta}$ and the fractional integral operator $I_{\gamma}$ (as defined by relation (\ref{riezs})) involving the spaces $\left( L^{q},L^{p}\right) ^{\alpha }(X)$ and weights fulfilling condition of $\mathcal A_{\infty}$ type stated in terms of Orlicz norms as in \cite{1}. 

 The main result is Theorem \ref{TH: cont op max} which is an extension of Theorem \ref{theoB} and contains, as a special case, the following result.
\begin{thm}\label{theoD}
Assume that\\
$\bullet$ there is a positive non decreasing function $\varphi$ defined on $\left]0,\infty\right[$ and positive constants $\mathfrak a$ and $\mathfrak b$ such that
\begin{equation}
\mathfrak a\varphi(r)\leq\mu\left(B_{(x,r)}\right)\leq \mathfrak b\varphi(r)\ \ x\in X\text{ and } 0<r,\label{normal}
\end{equation}
 $\bullet\;q,\alpha,p$ and $\beta$ are elements of $\left[1,\infty\right]$ such that $q\leq\alpha\leq p$ and $0<\frac{1}{s}=\frac{1}{\alpha}-
\frac{1}{\beta} \leq\frac{1}{q}-\frac{1}{\beta}\leq\frac{1}{p}$.

Then there is a real constant $C$ such that, for any $\mu$-measurable function $f$ on $X$ we have
\begin{equation}
\left\|\mathcal M_{q,\beta}f\right\|^{\ast}_{s,\infty}\equiv\sup_{\theta>0}\theta\mu\left(\left\{x\in X:\mathcal M_{q,\beta}f(x)>\theta\right\}\right)^{\frac{1}{s}}\leq C\left\|f\right\|_{q,p,\alpha}.\label{comp}
\end{equation}
\end{thm}
Remark that condition (\ref{normal}) is satisfied in the following cases:

$\bullet$ $X$ is an Ahlfors $n$ regular metric space, i. e. there is a positive integer $n$ and a positive constant $C$ which is independent of the main parameters such that $C^{-1}r^{n}\leq\mu\left(B_{(x,r)}\right)\leq C r^{n}$,

$\bullet$ $X$ is a Lie group with polynomial growth equipped with a left  Haar measure $\mu$ and the Carnot-Carath\'eodory metric $d$ associated with a H\"ormander system of left invariant vector fields (see \cite{HMY},\cite{Ma} and \cite{Va}).

Let us assume the hypotheses of Theorem \ref{theoD} and that $q<\alpha<p$. Theorem 2.12 of \cite{f-f} assert that weak-$L^{\alpha}(X)$ is strictly included in $(L^{q},L^{p})^{\alpha}(X)$. So, we may find an element $f_{0}$ in $(L^{q},L^{p})^{\alpha}(X)$ which is not in weak-$L^{\alpha}(X)$ space. Considering this element, Theorem \ref{theoD} asserts that $\mathcal M_{q,\beta}f_{0}$ belongs to the weak-$L^{s}$ space, while Theorem 2-7 of \cite{1} gives no control on it. This remark shows that, even if $\mathcal M_{q,\beta}$ is a particular case of the maximal operator $\mathcal M_{\psi}$ under consideration in Theorem 2-7 of \cite{1}, the range of application of this last theorem is different from that of our Theorem \ref{TH: cont op max}.

It is worth noting that $\mathcal M_{q,\beta}$ satisfies a norm inequality similar to (\ref{inC2}) (see Theorem \ref{th: cont rec op maximal}). This implies that if the maximal function $\mathcal M_{q,\beta}f$ belongs to some weak-Lebesgue space, then $f$ is in some $(L^{q},L^{p})^{\alpha}(X)$.

Let us consider the following fractional operator $I_{\gamma}$ $(0<\gamma<1)$ defined by 
\begin{equation}
I_{\gamma}f(x)=\int_{X}\frac{f(y)d\mu(y)}{\mu\left(B(x,d(x,y))\right)^{1-\gamma}}.\label{riezs}
\end{equation}
This operator is clearly an extension of the classical Riesz potential operator in $\mathbb R^{n}$. As in the euclidean case, $I_{\gamma}$ is controlled in norm by $\mathcal M_{1,\beta}$ where $\beta=\frac{1}{\gamma}$ (see Theorem \ref{TH: contrôle faib Intfract<opmax}). Thus from the weight norm inequality on $\mathcal M_{1,\beta}$ stated in Theorem \ref{TH: cont op max}, we may deduce a similar one on $I_{\gamma}$.

The remaining of the paper is organised as follows: section 2 is devoted to continuity properties of $\mathcal M_{q,\beta}$ and contains also background elements on homogeneous spaces, Young functions and  $(L^{q},L^{p})^{\alpha}(X)$ spaces. In section 3 we extend the results on $\mathcal M_{q,\beta}$ to $I_{\gamma}$. Throughout the paper, we will denote by $C$ a positive constant which is independent of the main parameters, but it may vary from line to line. Constant with subscripts, such as $C_{\mu}$ do not change in different occurences.

\section{Continuity of the fractional maximal operators $\mathcal M_{q,\protect\beta}$}
Let $X=(X,d,\mu)$ be a space of homogeneous type: $(X,d)$ is a quasi-metric space endowed with a non negative Borel measure $\mu $ satisfying the following doubling condition 
\begin{equation}
\mu \left( B_{\left( x,2r\right) }\right) \leq C\mu \left( B_{\left(
x,r\right) }\right) <\infty, \ \ x\in X\text{ and } r>0\label{0.001}
\end{equation}
 where $B_{\left( x,r\right) }=\left\{ y\in X:d(x,y)<r\right\}$ is the ball of center $x$ and radius $r$ in $X$. If $B$ is  an arbitrary ball, then  we denote by $x_{B}$ its center and $r(B)$ its radius, and for any real number $\delta>0$, $\delta B$ denotes the ball centered at $x_{B}$ with radius $\delta r(B)$.
 
  Since $d$ is a quasimetric, there exists a constant $\kappa\geq 1$ such that 
\begin{equation}
d(x,z)\leq\kappa\left(d(x,y)+d(y,z)\right),\ \ \ \ \ \ x,y,z\in X.
\end{equation}

If $C_{\mu}$ is the smallest constant for which (\ref{0.001}) 
holds, then $D_{\mu}=\log _{2}C_{\mu}$ is called the
doubling order of $\mu$. It is known (\cite{6} or \cite{4}) that for all balls  $B_{2}\subset B_{1}$ of $\left( X,d\right)$
\begin{equation}
\frac{\mu \left( B_{1}\right) }{\mu \left( B_{2}\right) }\leq C_{\mu}\left( \frac{r\left(
B_{1}\right) }{r\left( B_{2}\right) }\right) ^{D_{\mu}}.
\label{0.05}
\end{equation}
A quasimetric $\delta$ on $X$ is said to be equivalent to $d$ if there exist constants $C_{1}>0$ and $C_{2}>0$ such that
\begin{equation*}
C_{1}d(x,y)\leq\delta(x,y)\leq C_{2}d(x,y),\ \ x,y\in X.
\end{equation*}
We observe that topologies defined by equivalent quasi-metrics on $X$ are equivalent. It is shown in \cite{MS}, that there is a quasi-metric $\delta$ equivalent to $d$ for which balls are open sets.
\medskip

In the sequel we assume that $X=\left( X,d,\mu \right) $ is a fixed space of homogeneous type and\\
$\bullet$ all balls $B_{\left( x,r\right) }=\left\{y\in X:d\left( x,y\right) <r\right\} $ are open subsets of $X$ endowed with the $d$-topology and $(X,d)$ is separable,\\ 
$\bullet \ \mu(X)=\infty$\\
$\bullet \ \mu(\left\{x\right\})=0,\ x\in X,$\\
$\bullet \  B_{(x,R)}\setminus B_{(x,r)}\neq \emptyset,\ 0<r<R<\infty$, and $x\in X$, so that as proved in \cite{W}, there exist two constants $\tilde{C}_{\mu}>0$ and $\delta_{\mu}>0$ such that 
 \begin{equation}
 \frac{\mu(B_{1})}{\mu(B_{2})}\geq \tilde{C}_{\mu}\left(\frac{r(B_{1})}{r(B_{2})}\right)^{\delta_{\mu}}\text{ for all balls }B_{2}\subset B_{1}\text{ of }X.\label{revd}
 \end{equation}
We will now recall the concepts necessary to express the conditions we impose on our weights.

\begin{defn} 
Let $\Phi$ be a non negative function on $\left[ 0,\infty \right)$
\begin{enumerate}
\item[a)] $\Phi$ is a Young function if it is continuous, non decreasing, convex and satisfies the conditions 
\begin{equation*}
\Phi (0)=0\ \text{and\ }\underset{x\rightarrow \infty }{\lim }\Phi
(x)=\infty.
\end{equation*}
\item[b)] Assume that $\Phi$ is a Young function:
\begin{enumerate}
\item[(i)]it is doubling if there is a constant $C>0$ such that   
\begin{equation*}
\Phi \left( 2t\right) \leq C\Phi \left( t\right) \ \text{for all} \ t\geq0,
\end{equation*}
 
\item[(ii)]it satisfies the $B_{p}$ condition ($1\leq p<\infty$) if there is a number $a>0$ such that  
\begin{equation*}
\int^{\infty }_{a}\frac{\Phi (t)}{t^{p}}\frac{dt}{t}<\infty,
\end{equation*}
\item[(iii)] its conjugate $\Phi^{\ast }$, is defined by  
\begin{equation}
\Phi ^{\ast }(u)=\sup \left\{ tu-\Phi(t):t\in \mathbb{R}%
_{+}\right\},
\end{equation}
\item[(iv)] for any $\mu$-measurable function $f$ on $X$ 
\begin{equation}
\left\| f\right\| _{\Phi ,B}=\inf \left\{ a>0:\frac{1}{\mu \left( B\right) }%
\int_{B}\Phi \left( a^{-1}\left| f\right| \right) d\mu \leq 1\right\}
\end{equation}
for any ball $B$ in $X$, and
\begin{equation}
M_{\Phi}f(x)=\sup_{ball\ B\ni x}\left\| f\right\| _{\Phi ,B}.
\end{equation} 
\end{enumerate}
\end{enumerate}
\end{defn}
It is proved in Theorem 5.1 of \cite{2} that a doubling Young function $\Phi$ belongs to the class $B_{p}$ with $1<p<\infty$ if and only if there exists a constant $C>0$ such that
\begin{equation}
\int_{X}\left(\mathcal M_{\Phi}f(x)\right)^{p}d\mu(x)\leq C\int_{X}f(x)^{p}d\mu(x)\label{27}
\end{equation}
for all nonnegative $f$. We also have the local version of generalized H$\ddot{\text{o}}$lder inequality 
 
\begin{equation}
\frac{1}{\mu \left( B\right) }\int_{B}\left| fg\right| d\mu \leq \left\|
f\right\| _{\Phi ,B}\left\| g\right\| _{\Phi ^{\ast },B},  \label{7}
\end{equation}
which is valid for all measurable functions $f$ and $g$, and for all ball $B$. For more information about Young function, see \cite{RR}.

We will need the following covering lemma stated and proved in \cite{6}.

\begin{lem}
\label{lem:famiden}Let $\mathfrak{F}$\ be a family of balls with bounded
radii. Then, there exists a countable subfamily of disjoint balls $\left\{
B_{\left( x_{i},r_{i}\right) },i\in J\right\} $ such that each ball in $%
\mathfrak{F}$ is contained in one of the balls $B_{\left( x_{i},3\kappa
^{2}r_{i}\right) }$, for some $i\in J$.
\end{lem}
We are now ready to state and prove our main result.
\begin{thm}
\label{TH: cont op max}\ Let $ q,\alpha ,p,q_{1},\alpha
_{1},p_{1,}\beta  $ be elements of $\left[ 1,\infty \right]
$ such that 
\begin{equation*}
1\leq q\leq \alpha \leq p\ \text{with\ }0<\frac{1}{\alpha }-\frac{1}{%
\beta }=\frac{1}{s},
\end{equation*}
and 
\begin{equation*}
q<q_{1}\leq \alpha _{1}\leq p_{1}<\infty \ \text{with\ }0<\frac{1}{q_{1}}-%
\frac{1}{\beta }=\frac{1}{t}\leq \frac{1}{p_{1}}.
\end{equation*}
Let $\left( w,v\right) $ be a pair of weights for which there exists a
constant $A$ such that 
\begin{equation}
\mu \left( B\right) ^{-\frac{1}{t}}\left\| w\chi _{B}\right\|
_{t}\left\| v^{-q}\right\| _{\Phi ,B}^{\frac{1}{q}}
\leq A  \label{1}
\end{equation}
for all balls B in $\left( X,d\right) $, where $\Phi $ is a doubling Young
function whose conjugate function $\Phi ^{\ast }$ satisfies the $B_{\frac{%
q_{1}}{q}}$ condition. Then there is a constant $C$ such that for all
$\mu$-measurable function $f$, and $\theta >0$, we have 
\begin{equation}
\left(\int_{\Pi _{\theta }}w^{t}(x)d\mu(x)\right)^{\frac{1}{t}}\leq C\theta^{-1}\left\|fv\right\|_{q_{1}}, \label{wheeden}
\end{equation}
and if we assume that $\mu$ satisfies condition (\ref{normal}), then 
\begin{equation}
\left(\int_{\Pi _{\theta }}w^{t}(x)d\mu (x)\right) ^{\frac{1}{t}}\leq
C\left( \theta ^{-1}\left\| fv\right\| _{q_{1},p_{1},\alpha _{1}}\right)
\left( \theta ^{-1}\left\| f\right\| _{q,p,\alpha }\right) ^{s\left( \frac{1%
}{q_{1}}-\frac{1}{\alpha _{1}}\right) },\label{nouveaux}
\end{equation}
where $\Pi _{\theta }=\left\{ x\in X:\mathcal M_{q,\beta
}f(x)>\theta \right\} .$
\end{thm}

\proof
Inequality (\ref{wheeden}) is immediate from Theorem 2.7 of \cite{1}. We just have to prove inequality (\ref{nouveaux}).
  
Let $f$ be an element of $\left( L^{q},L^{p}\right) ^{\alpha }(X)$. Fix  $\theta >0$. For $x$ in $\Pi _{\theta }$, there exists $r_{x}$ such that
\begin{equation}
\mu\left(B_{(x,r_{x})}\right) ^{\frac{1}{\beta }-\frac{1}{q}
}\left\| f\chi _{_{B_{\left( x,r_{x}\right) }}}\right\| _{q}>\theta 
\label{I1}
\end{equation}
and therefore
\begin{equation}
\mu\left(B_{(x,r_{x}}\right)\leq\left(\theta^{-1}\left\|f\right\|_{q,\infty,\alpha}\right)^{s}.\label{I2}
\end{equation}
Fix a ball $B_{(x_{0},R)}$ in $X$ and set $\Pi^{R}_{\theta}=\Pi_{\theta}\cap B_{(x_{0},R)}$. For any $x$ in $\Pi^{R}_{\theta}$ we have $B_{(x_{0},R)}\subset B_{(x,r_{x})}$ provided $r_{x}>2\kappa R$. It follows from the reverse doubling property and (\ref{I2}) that 
\begin{equation}
r^{\delta_{\mu}}_{x}\leq C^{-1}_{\mu}\frac{R^{\delta_{\mu}}}{\mu\left(B_{(x_{0},R)}\right)}\left(\theta^{-1}\left\|f\right\|_{q,\infty,\alpha}\right)^{s}.\label{I3}
\end{equation}
So we obtain that for any $x$ in $\Pi^{R}_{\theta}$ 
\begin{equation}
r^{\delta_{\mu}}_{x}\leq\max\left\{2\kappa R, C^{-1}_{\mu}\frac{R^{\delta_{\mu}}}{\mu\left(B_{(x_{0},R)}\right)}\left(\theta^{-1}\left\|f\right\|_{q,\infty,\alpha}\right)^{s}\right\}<\infty.
\end{equation}
Thus by Lemma \ref{lem:famiden}, the family $\mathcal F=\left\{B_{(x,r_{x})}:x\in\Pi^{R}_{\theta}\right\}$ has a countable subfamily $\left\{B_{i}:i\in J\right\}$ of disjoint balls such that each element $B$ of $\mathcal F$ is contained in some $3\kappa^{2}B_{i}$.

Let $i$ be an element of $J$. By (\ref{I1}) and the generalized H\"older inequality we have
\begin{eqnarray*}
\theta ^{q}&\leq&\mu \left( B_{i}\right) ^{\frac{q}{\beta }}\left(\frac{1}{\mu\left(B_{i}\right)}\int_{B_{i}}\left|fvv^{-1}\right|d\mu\right)\\
&\leq& C\mu \left( B_{i}\right) ^{\frac{q}{\beta }}\left\|(fv\chi_{B_{i}})^{q}\right\|_{\Phi^{\ast},3\kappa^{2}B_{i}}\left\|v^{-1}\chi_{B_{i}}\right\|_{\Phi,3\kappa^{2}B_{i}}\\
&\leq&C\mu \left( B_{i}\right) ^{\frac{q}{\beta }}M_{\Phi ^{\ast }}\left( fv\chi _{B_{i}}\right) ^{q}(y)\left\| v^{-q}\chi _{B_{i}}\right\| _{\Phi ,3\kappa ^{2}B_{i}}\\
\end{eqnarray*}
for any $y$ in $B_{i}$. So we obtain 
\begin{equation}
\theta ^{q}\mu \left( B_{i}\right)\leq C\mu \left( B_{i}\right) ^{\frac{q}{\beta }}\int_{B_{i}}M_{\Phi ^{\ast }}\left( fv\chi _{B_{i}}\right) ^{q}(y)d\mu(y)\left\| v^{-q}\chi _{B_{i}}\right\| _{\Phi ,3\kappa ^{2}B_{i}}.
\end{equation}
 Applying H\"older inequality and (\ref{27}) we get
 \begin{eqnarray*}
 \theta ^{q}&\leq&C\mu \left( B_{i}\right) ^{-\frac{q}{t }}\left[\int_{B_{i}}\left\{M_{\Phi ^{\ast }}\left( fv\chi _{B_{i}}\right) ^{q}(y)\right\}^{\frac{q_{1}}{q}}d\mu(y)\right]^{\frac{q}{q_{1}}}\left\| v^{-q}\chi _{B_{i}}\right\| _{\Phi ,3\kappa ^{2}B_{i}}\\
 &\leq&C\mu \left( B_{i}\right) ^{-\frac{q}{t }}\left\|fv\chi _{B_{i}}\right\|^{q}_{q_{1}}\left\| v^{-q}\chi _{B_{i}}\right\| _{\Phi ,3\kappa ^{2}B_{i}};
 \end{eqnarray*}
 that is
 \begin{equation}
 1\leq C\theta^{-1}\mu \left( B_{i}\right) ^{-\frac{1}{t }}\left\|fv\chi _{B_{i}}\right\|_{q_{1}}\left\| v^{-q}\chi _{B_{i}}\right\|^{\frac{1}{q}} _{\Phi ,3\kappa ^{2}B_{i}}.\label{I4}
 \end{equation}
As $\Pi^{R}_{\theta}\subset\cup_{i\in J}3\kappa^{2}B_{i}$ and $\frac{p_{1}}{t}\leq 1$, we have 
\begin{eqnarray*}
\left\|w\chi_{\Pi^{R}_{\theta}}\right\|_{t}&\leq&\left(\sum_{i\in J}\left\|w\chi_{3\kappa^{2}B_{i}}\right\|^{p_{1}}_{t}\right)^{\frac{1}{p_{1}}}\\
&\leq&C\theta^{-1}\left[\sum_{i\in J}\left(\mu\left(B_{i}\right)^{-\frac{1}{t}}\left\|fv\chi_{B_{i}}\right\|_{q_{1}}\left\|v^{-q}\chi_{B_{i}}\right\|^{\frac{1}{q}}_{\Phi,3\kappa^{2}B_{i}}\left\|w\chi_{3\kappa^{2}B_{i}}\right\|_{t}\right)^{p_{1}}\right]^{\frac{1}{p_{1}}}.
\end{eqnarray*}
Thus, according to assumption (\ref{I2})
\begin{equation}
\left\|w\chi_{\Pi^{R}_{\theta}}\right\|_{t}\leq C\theta^{-1}\left(\sum_{i\in J}\left\|fv\chi_{B_{i}}\right\|^{p_{1}}_{q_{1}}\right)^{\frac{1}{p_{1}}}.\label{*}
\end{equation}

Let $n$ be a positive integer and set

$\bullet\ J_{n}=\left\{i\in J:\frac{1}{n}\leq r(B_{i})\right\}$

$\bullet\ m_{n}$ and $\bar{k}$ the integers satisfying respectively 
\begin{equation}
\rho^{m_{n}+1}\leq\frac{1}{2\kappa n}<\rho^{m_{n}+2}\text{ and } \rho^{\bar{k}+1}\leq\frac{r}{2\kappa}<\rho^{\bar{k}+2},
\end{equation}
where $r=\sup\left\{r(B_{i}), i\in I\right\}$ and $\rho=8\kappa^{5}$.

It is proved in \cite{sw} that there are points $x^{k}_{j}$ and Borel sets $E^{k}_{j}$, $1\leq j<N_{k}$, $k\geq m_{n}$ (where $N_{k}\in\mathbb N\cup\left\{\infty\right\}$), such that 
\begin{enumerate}
\item[(i)]$B_{\left( x^{k}_{j},\rho ^{k}\right) }\subset E^{k}_{j}\subset
B_{\left( x^{k}_{j},\rho ^{k+1}\right) },\ 1\leq j<N_{k},\;k\geq m_{n}$
\item[(ii)] $X=\cup_{j}E^{k}_{j} \;k\geq m_{n}$ and $E^{k}_{j}\cap E^{k}_{i}=\emptyset$ if $i\neq j$
\item [(iii)]given $i,j,k,\ell$ with $m_{n}\leq k<\ell$, then either $E^{k}_{j}\subset E^{\ell}_{i}$ or $E^{k}_{j}\cap E^{\ell}_{i}=\emptyset$.
\end{enumerate}
Let $i$ be an element of $I_{n}$. Denote by $k_{i}$ the integer satisfying 
\begin{equation}
\rho^{k_{i}+1}\leq\frac{r(B_{i})}{2\kappa}<\rho^{k_{i}+2}
\end{equation}
 and set $ L_{i}=\left\{j:1\leq j<N_{k_{i}},\;E^{k_{i}}_{j}\cap B_{i}\neq\emptyset\right\}$.

We know that the number of elements of $L_{i}$ is less than a constant $\mathfrak N$ depending only on the structure constants $(\kappa,C_{\mu},D_{\mu},\tilde{C}_{\mu},\delta_{\mu})$ (see inequality (43) of \cite{f-f}).
Denoting by $j_{i}$ an element of $L_{i}$ satisfying 
\begin{equation}
\left\|fv\chi_{B_{i}\cap E^{k}_{j_{i}}}\right\|_{q_{1}}=\max_{j\in L_{i}}\left\|fv\chi_{B_{i}\cap E^{k}_{j}}\right\|_{q_{1}},
\end{equation}
we have
\begin{equation}
\left\|fv\chi_{B_{i}}\right\|_{q_{1}}\leq\mathfrak N^{-1}\left\|fv\chi_{B_{i}\cap E^{k}_{j_{i}}}\right\|_{q_{1}}.
\end{equation}
Hence
\begin{eqnarray*}
\left(\sum_{i\in J_{n}}\left\|fv\chi_{B_{i}}\right\|^{p_{1}}_{q_{1}}\right)^{\frac{1}{p_{1}}}&\leq&\mathfrak N^{-1}\left(\sum_{i\in J_{n}}\left\|fv\chi_{E^{k_{i}}_{j_{i}}\cap B_{i}}\right\|^{p_{1}}_{q_{1}}\right)^{\frac{1}{p_{1}}}\\
&=&\mathfrak N^{-1}\left(\sum^{N_{\bar{k}}}_{\ell=1}\sum_{ i\in J_{n}: 
E^{k_{i}}_{j_{i}}\subset E^{\bar{k}}_{\ell}}\left\|fv\chi_{E^{k_{i}}_{j_{i}}\cap B_{i}}\right\|^{p_{1}}_{q_{1}}\right)^{\frac{1}{p_{1}}}\\
&\leq&\mathfrak N^{-1}\left(\sum^{N_{\bar{k}}}_{\ell=1}\left\|fv\chi_{E^{\bar{k}}_{j_{i}}\cap (\cup_{i\in J_{n}}B_{i}}\right\|^{p_{1}}_{q_{1}}\right)^{\frac{1}{p_{1}}}\\
&=&\mathfrak N^{-1}\left[\sum^{N_{\bar{k}}}_{\ell=1}\left(\mu\left(E^{\bar{k}}_{\ell}\right)^{\frac{1}{\alpha_{1}}-\frac{1}{q_{1}}}\left\|fv\chi_{E^{\bar{k}}_{j_{i}}\cap (\cup_{i\in J_{n}}B_{i})}\right\|_{q_{1}}\right)^{p_{1}}\mu\left(E^{\bar{k}}_{\ell}\right)^{\frac{p}{q_{1}}-\frac{p}{\alpha_{1}}}\right]^{\frac{1}{p_{1}}}.
\end{eqnarray*}
Notice that for any $1\leq \ell<N_{\bar{k}}$, we have
\begin{eqnarray*}
\mu\left(E^{\bar{k}}_{\ell}\right)&\leq&\mu\left(B_{(x^{\bar{k}}_{\ell},\rho^{\bar{k}+1})}\right)\leq\mathfrak b\varphi(\rho^{\bar{k}+1})\\
&\leq&\mathfrak b\varphi\left(\frac{r}{2\kappa}\right)\leq\mathfrak b\sup_{i\in J}\varphi\left(\frac{r(B_{i})}{2\kappa}\right)\\
&\leq&\mathfrak b\sup_{i\in J}\mathfrak a^{-1}\mu\left(\frac{1}{2\kappa}B_{i}\right)\leq \mathfrak b\mathfrak a^{-1}\left(\theta^{-1}\left\|f\right\|_{q,\infty,\alpha}\right)^{s}.
\end{eqnarray*}
Therefore,
\begin{equation}
\left(\sum_{i\in J_{n}}\left\|fv\chi_{B_{i}}\right\|^{p_{1}}_{q_{1}}\right)^{\frac{1}{p_{1}}}\leq C\left[\sum^{N_{\bar{k}}}_{\ell=1}\left(\mu\left(E^{\bar{k}}_{\ell}\right)^{\frac{1}{\alpha_{1}}-\frac{1}{q_{1}}}\left\|fv\chi_{E^{\bar{k}}_{j_{i}}}\right\|_{q_{1}}\right)^{p_{1}}\right]^{\frac{1}{p_{1}}}\left(\theta^{-1}\left\|f\right\|_{q,\infty,\alpha}\right)^{s\left(\frac{1}{q_{1}}-\frac{1}{\alpha_{1}}\right)}.\label{39}
\end{equation}
Since the last formula does not depend on $n$, we get from (\ref{*}) and (\ref{39})
\begin{equation}
\left\|w\chi_{\Pi^{R}_{\theta}}\right\|_{t}\leq\theta^{-1}\left[\sum^{N_{\bar{k}}}_{\ell=1}\left(\mu\left(E^{\bar{k}}_{\ell}\right)^{\frac{1}{\alpha_{1}}-\frac{1}{q_{1}}}\left\|fv\chi_{E^{\bar{k}}_{j_{i}}}\right\|_{q_{1}}\right)^{p_{1}}\right]^{\frac{1}{p_{1}}}\left(\theta^{-1}\left\|f\right\|_{q,\infty,\alpha}\right)^{s\left(\frac{1}{q_{1}}-\frac{1}{\alpha_{1}}\right)}.
\end{equation}
We recall that Proposition 4.1 of \cite{f-f} asserts that there are positive constants $C_{1}$ and $C_{2}$ not depending on $r$ and $fv$ such that 
\begin{equation}
C_{1}\;_{r}\left\|fv\right\|_{q_{1},p_{1},\alpha_{1}}\leq\left[\sum^{N_{\bar{k}}}_{\ell=1}\left(\mu\left(E^{\bar{k}}_{\ell}\right)^{\frac{1}{\alpha_{1}}-\frac{1}{q_{1}}}\left\|fv\chi_{E^{\bar{k}}_{j_{i}}}\right\|_{q_{1}}\right)^{p_{1}}\right]^{\frac{1}{p_{1}}}\leq C_{2}\;_{r}\left\|fv\right\|_{q_{1},p_{1},\alpha_{1}}.
\end{equation}
So we have
\begin{eqnarray*}
\left\|w\chi_{\Pi^{R}_{\theta}}\right\|_{t}&\leq& C\theta^{-1}\;_{r}\left\|fv\right\|_{q_{1},p_{1},\alpha_{1}}\left(\theta^{-1}\left\|f\right\|_{q,\infty,\alpha}\right)^{s\left(\frac{1}{q_{1}}-\frac{1}{\alpha_{1}}\right)}\\
&\leq& C\theta^{-1}\left\|fv\right\|_{q_{1},p_{1},\alpha_{1}}\left(\theta^{-1}\left\|f\right\|_{q,\infty,\alpha}\right)^{s\left(\frac{1}{q_{1}}-\frac{1}{\alpha_{1}}\right)}.
\end{eqnarray*}
As $(x_{0},R)$ is arbitrary in $X\times\left]0,\infty\right[$, we obtain 
\begin{equation}
\left\|w\chi_{\Pi_{\theta}}\right\|_{t}\leq C\theta^{-1}\left\|fv\right\|_{q_{1},p_{1},\alpha_{1}}\left(\theta^{-1}\left\|f\right\|_{q,\infty,\alpha}\right)^{s\left(\frac{1}{q_{1}}-\frac{1}{\alpha_{1}}\right)}.
\end{equation}
\epf

In the proof of the above theorem, the condition $q<q_{1}$ is needed only when we have to use the $B_{\frac{q_{1}}{q}}$ characterization. When $w=v=1$ this characterization is not needed. So Theorem \ref{theoD} follows immediately  from Theorem \ref{TH: cont op max}. 

The next theorem is some kind of reverse for Theorem\ref{theoD}

\begin{thm}
\label{th: cont rec op maximal}Let $q,\alpha $,$u$ and $v$ be elements of $%
\left[ 1,\infty \right] $ such that  
\begin{equation*}
q\leq \alpha ,\ 0\leq \frac{1}{\alpha }-\frac{1}{\beta }=\frac{1}{s}\ 
\text{and }u\leq s\leq v.
\end{equation*}
Then there is a constant $D$ such that for any $\mu$-measurable function $f$ 
\begin{equation}
\left\| f\right\| _{q,v,\alpha }\leq D\left\| \mathcal M_{q,\beta }f\right\| _{u,v,s}.
\label{1.1}
\end{equation}
\end{thm}

\proof
 Let $f$ be such that $\left\|\mathcal M_{q,\beta}f\right\|_{u,v,s}<\infty$. We notice that under the hypothesis, we have $q\leq \alpha\leq s\leq v$ and $\alpha\leq\beta$.
 
 $1^{rst}$ case: $q=\infty$. Then $\alpha =\beta =s=v=\infty $ and therefore, it follows from the definitions that
\begin{equation*}
\left\| f\right\| _{\infty ,\infty ,\infty }=\left\| f\right\| _{\infty
}=\left\|\mathcal M_{\infty ,\infty }f\right\| _{\infty }\leq C\left\| \mathcal M_{\infty
,\infty }f\right\| _{u,\infty ,\infty }.
\end{equation*}

$2^{nd}$ case: $q<\infty$.
\begin{enumerate}
\item[(a)]If $u=\infty $, then $s=v=\infty $, $\alpha =\beta $ and therefore, 
\begin{eqnarray*}
\left\| f\right\| _{q,\infty ,\alpha }&=&\sup_{r>0}\sup_{x\in X}\text{ess} \;\mu\left(B_{(x,r)}\right)^{\frac{1}{\alpha }-\frac{1}{q}}\left\| f\chi _{_{B_{\left( x,r\right) }}}\right\| _{q}\\
&=&\sup_{r>0}\sup_{x\in X}\text{ess} \;\mu\left(B_{(x,r)}\right)^{\frac{1}{\beta }-\frac{1}{q}}\left\| f\chi _{_{B_{\left( x,r\right) }}}\right\| _{q}=\left\|\mathcal M_{q,\beta }f\right\|
_{\infty}=\left\|\mathcal M_{q,\beta }f\right\|
_{\infty ,\infty ,\infty }.
\end{eqnarray*}
\item[(b)] Suppose that $u<\infty ,$ and consider two positive reals numbers  $r$ and $r_{_{1}}$  satisfying $r_{_{1}}=\dfrac{r}{2\kappa }$. For any $y\in X$ and $x\in B_{(y,r_{1})}$, we have $B_{(y,r_{1})}\subset B_{(x,r)}$ and therefore, by the doubling condition
\begin{equation}
\mathcal M_{q,\beta}f(x)\geq C_{\mu}^{\frac{1}{\beta}-\frac{1}{q}}\mu\left(B_{(y,r_{1})}\right)^{\frac{1}{\beta}-\frac{1}{q}}\left\|f\chi_{_{B_{(y,r_{1})}}}\right\|_{q}.
\end{equation}
From this, it follows that for any $y\in X$, we have
\begin{equation}
\left\|\mathcal M_{q,\beta}f\chi_{B_{(y,r_{1})}}\right\|_{u}\geq C_{\mu}^{\frac{1}{\beta}-\frac{1}{q}}\mu\left(B_{(y,r_{1})}\right)^{\frac{1}{u}+\frac{1}{\beta}-\frac{1}{q}}\left\|f\chi_{_{B_{(y,r_{1})}}}\right\|_{q}
\end{equation}
and therefore, 
\begin{equation}
\mu\left(B_{(y,r_{1})}\right)^{\frac{1}{s}-\frac{1}{v}-\frac{1}{u}}\left\|\mathcal M_{q,\beta}f\chi_{B_{(y,r_{1})}}\right\|_{u}\geq C_{\mu}^{\frac{1}{\beta}-\frac{1}{q}}\mu\left(B_{(y,r_{1})}\right)^{\frac{1}{\alpha}-\frac{1}{v}-\frac{1}{q}}\left\|f\chi_{_{B_{(y,r_{1})}}}\right\|_{q}.
\end{equation}
This yields immediately the desired inequality
\end{enumerate}
\epf

\section{Continuity of the fractional integral $I_{\protect\alpha }f$}
It is known in the euclidean case that the fractional integral $I_{\gamma}f$ is controlled in norm by the fractional maximal function $\mathfrak m_{1,\frac{1}{\gamma}}f$ (seeTheorem 1 of \cite{3}). We give the analogous of this control in the setting of spaces of homogeneous type.

\begin{thm}
\label{TH: contrôle faib Intfract<opmax}Let $0<q<\infty $, $0<\gamma <1$
and a weight $w$ in $\mathcal{A}_{\infty }.$ There is a constant $C$ such
that for any $\mu$-measurable function $f$  
\begin{equation*}
\underset{a>0}{\sup }a^{q}\int_{E_{a}}w(x)d\mu (x)\leq C\underset{a>0}{\sup 
}a^{q}\int_{F_{a}}w(x)d\mu (x),
\end{equation*}
where$\ E_{a}=\left\{ x\in X:\left| I_{\gamma }f(x)\right|
>a\right\} $ and $F_{a}=\left\{ x\in X:\mathcal M_{1,\frac{1}{\gamma }%
}f(x)>a\right\} .$
\end{thm}

\proof
In our argumentation, we shall adapt the proof of Theorem 1 of \cite{3}, keeping in mind that we do not have a Withney decomposition avalaible.
\begin{enumerate}
\item[1)] Let $f$ be a $\mu$-measurable, non negative, bounded function, with a support included in a ball $B_{0}=B_{(x_{0},k_{0})}$. According to Lemma 6 of \cite{4}, there exists a constant $C_{0}>0$ not depending on $f$, such that 
\begin{equation*}
 I_{\gamma}f\leq \mathcal M\left( I_{\gamma}f\right) \leq C_{0}I_{\gamma }f,
\end{equation*}
where $\mathcal M=\mathcal M_{1,\infty}$. Let $\theta$ be a positive number and set 
\begin{equation*}
\widetilde{E}_{\theta }=\left\{ x\in X: \mathcal M(I_{\gamma}f)(x) >\theta \right\} \text{ and } E_{\theta }=\left\{ x\in X: I_{\gamma}f(x) >\theta\right\}.
\end{equation*}
The set $E_{\theta }$ is included in $\widetilde{E}_{\theta }$
which is opened and satisfy $\mu \left( \widetilde{E}_{\theta
}\right) <\infty .$ According to Lemma 8 of \cite{4}, there exists a countable family $\left\{ B_{\left( x_{i},r_{i}\right) };i\in J\right\} $  of pairwise disjoint balls and two positive constants $M$ and $c$ depending only on the structure constants of $X$, such that
\begin{equation}
\left\{ 
\begin{array}{l}
\widetilde{E}_{\theta }=\cup _{i\in J}B_{\left( x_{i},cr_{i}\right) }, \\ 
\underset{i\in J}{\sum }\chi _{B_{\left( x_{i},2\kappa cr_{i}\right) }}\leq
M\chi _{\widetilde{E}_{\theta }}, \\ 
B_{\left( x_{i},4\kappa ^{2}cr_{i}\right) }\cap \left( X\setminus 
\widetilde{E}_{\theta }\right) \neq \emptyset \ \text{for all}\  i\in J.%
\end{array}
\right. \label{R1}
\end{equation}

Let us consider an element $\left( a,\varepsilon \right) $ of  $\left] 1,\infty \right[ \times \left] 0,1\right] $, and set 
$F_{\theta ,\varepsilon }=\left\{ x\in X:\mathcal M_{1,\frac{1}{\gamma}}f(x)>\theta \varepsilon \right\}$, $J_{1}=\left\{ i\in J:B_{\left( x_{i},cr_{i}\right) }\subset
F_{\theta ,\varepsilon }\right\}$ and $J_{2}=I\setminus J_{1}=\left\{ i\in J:B_{\left(x_{i},cr_{i}\right) }\setminus F_{\theta ,\varepsilon }\neq \emptyset
\right\}.$

Arguing as in the proof of Lemma 1 of \cite{3}, we obtain two constants $K>0$ and $B>1$ depending only on the structure constants of $X$, such that if $a\geq B$ and $i\in J_{2}$ then
\begin{equation}
\mu\left(\left\{x\in B_{(x_{i},cr_{i})}:I_{\gamma}f(x)>a\theta\right\}\right)\leq K\mu\left(B_{(x_{i},cr_{i})}\right)\left(\frac{\epsilon}{a}\right)^{\frac{1}{1-\gamma}}\label{R2}.
\end{equation}
Since 
\begin{equation*}
E_{\theta a}\subset E_{\theta }\subset \widetilde{E}_{\theta }=\underset{i\in J}{\cup }%
B_{\left( x_{i},cr_{i}\right) },
\end{equation*}
we have 
\begin{equation*}
E_{\theta a}=\left[ \underset{i\in J_{1}}{\cup }\left( E_{\theta a}\cap
B_{\left( x_{i},cr_{i}\right) }\right) \right] \cup \left[ \underset{i\in
J_{2}}{\cup }\left( E_{\theta a}\cap B_{\left( x_{i},cr_{i}\right) }\right)\right]\subset F_{\theta ,\varepsilon }\cup\left[\cup_{i\in J_{2}}\left(E_{\theta a\setminus F_{\theta\epsilon}}\right)\cap B_{(x_{i},cr_{i})}\right]
\end{equation*}
and therefore 
\begin{equation}
\int_{E_{\theta a}}w(x)d\mu (x)\leq  \int_{F_{\theta ,\varepsilon }}w(x)d\mu (x)+\underset{i\in J_{2}}{\sum }\int_{\left( E_{\theta a}\setminus
F_{\theta ,\varepsilon }\right)\cap B_{\left( x_{i},cr_{i}\right)}}w(x)d\mu (x).\label{R3}
\end{equation}

Now fix $a\geq B$ and $\rho >0$. Since $w$ is in $\mathcal A_{\infty}$, there exists $\delta>0$ such that for any ball $B$ in $X$ and any subset $E$ of $B$ satisfying $\mu \left( E\right) \leq \delta \mu \left( B\right)$, we have%
\begin{equation}
\int_{E}w(x)d\mu (x)\leq \rho \int_{B}w(x)d\mu (x).
\end{equation}
Choose $\overline{\varepsilon }\in \left] 0, 1\right]$ such that
$K\left( \frac{\overline{\varepsilon }}{a}\right) ^{\frac{1}{1-\gamma }%
}<\delta $ and take $0<\varepsilon <\min\left( \overline{\varepsilon },\frac{1}{C_{0}L}\right)$, where $L=C_{\mu}\left(2\kappa+4\kappa^{2}\right)^{1-\gamma}$. According to (\ref{R2}) we have for any  $i\in J_{2}$, 
\begin{equation}
\mu \left( B_{\left( x_{i},cr_{i}\right) }\cap E_{\theta
a}\right) <\delta \mu \left( B_{\left( x_{i},cr_{i}\right) }\right) 
\end{equation}
and therefore
\begin{equation}
 \int_{B_{\left( x_{i},cr_{i}\right) }\cap E_{\theta
a}}w(x)d\mu \left( x\right) \leq \rho \int_{B_{\left( x_{i},cr_{i}\right)
}}w(x)\mu \left( x\right) .
\end{equation}
From this inequality, (\ref{R3}) and (\ref{R1}) we obtain
 
\begin{equation}
\int_{E_{\theta a}}w(x)d\mu (x)\leq \int_{F_{\theta ,\varepsilon }}w(x)d\mu (x)+\rho M\int_{\widetilde{E}_{\theta }}w(x)d\mu (x).\label{R4}
\end{equation}
Let $x\in X\setminus 3\kappa B_{0}$. Assume that $0<t<\frac{1}{2}\inf_{y\in B_{0}}d(x,y)$ and $u_{t}\in B_{0}$ satisfies 
\begin{equation}
d(x,u_{t})-t\leq d(x,y),\ y\in B_{0}.
\end{equation}
We have 
\begin{equation}
2r(B_{0})\leq d(x,y)\leq\kappa\left[d(x,u_{t})+2\kappa r(B_{0})\right],\ y\in B_{0}
\end{equation}
and therefore
\begin{eqnarray*}
I_{\gamma}f(x)&\leq&\int_{B_{0}}\frac{f(y)}{\mu\left(B_{(x,d(x,y))}\right)^{1-\gamma}}d\mu(y)\leq\frac{1}{\mu\left(B_{(x,d(x,u_{t})-t)}\right)^{1-\gamma}}\int_{B_{(x,\kappa\left(d(x,u_{t})+2\kappa r(B_{0})\right))}}f(y)d\mu(y)\\
& \leq & C_{\mu}\left[\frac{\kappa\left(d(x,u_{t})+2\kappa r(B_{0})\right)}{d(x,u_{t})-t}\right]^{(1-\gamma)D_{\mu}}\mathcal M_{1,\frac{1}{\gamma}}f(x)\leq L\mathcal M_{1,\frac{1}{\gamma}}f(x).
\end{eqnarray*}
Hence,
\begin{equation}
\widetilde{E}_{\theta }\subset E_{\frac{\theta }{C_{0}}}\subset \left( E_{
\frac{\theta }{C_{0}}}\cap 3\kappa B_{0}\right) \cup F_{\theta ,\frac{1}{C_{0}L}}.
\end{equation}
We obtain from (\ref{R4})
\begin{eqnarray*}
\int_{E_{\theta a}}w(x)d\mu (x)&\leq& \int_{F_{\theta ,\varepsilon
}}w(x)d\mu (x)+\rho M\int_{E_{\frac{\theta }{C_{0}}}\cap 3\kappa
B_{0}}w(x)d\mu (x)+\rho M\int_{F_{\theta ,\frac{1}{C_{0}L}%
}}w(x)d\mu (x)\\
&\leq& \left( 1+\rho M\right) \int_{F_{\theta
,\varepsilon }}w(x)d\mu (x)+\rho M\int_{E_{\frac{\theta }{C_{0}}}\cap
3\kappa B_{0}}w(x)d\mu (x).
\end{eqnarray*}
That is  
\begin{eqnarray*}
\left( \theta a\right) ^{q}\int_{E_{\theta a}}w(x)d\mu (x) & \leq & \left(
1+\rho M\right) \left( \frac{a}{\varepsilon }\right) ^{q}\left( \theta
\varepsilon \right) ^{q}\int_{F_{\theta ,\varepsilon }}w(x)d\mu (x) \\ 
& + & \rho M\left( \frac{\theta }{C_{0}}\right) ^{q}\left( C_{0}a\right)
^{q}\int_{E_{\frac{\theta }{C_{0}}}\cap 3\kappa B_{0}}w(x)d\mu (x).
\end{eqnarray*}
Let $N$ be a positive integer. From the preceeding inequality we obtain
\begin{eqnarray*}
\underset{0<s<N}{\sup }s^{q}\int_{E_{s}}w(x)d\mu (x) & \leq & \left( 1+\rho
M\right) \left( \frac{a}{\varepsilon }\right) ^{q}\underset{0<s<N\frac{%
\varepsilon }{a}}{\sup }s^{q}\int_{F_{s,1}}w(x)d\mu (x) \\ 
& + & \rho M\left( C_{0}a\right) ^{q}\underset{0<s<\frac{N}{aC_{0}}}{\sup }%
s^{q}\int_{E_{s}\cap 3\kappa B_{0}}w(x)d\mu (x).%
\end{eqnarray*}
As 
$$\underset{0<s<\frac{N}{aC_{0}}}{\sup }s^{q}\int_{E_{s}\cap 3\kappa
B_{0}}w(x)d\mu (x)\leq \underset{0<s<N}{\sup }s^{q}\int_{E_{s}\cap 3\kappa
B_{0}}w(x)d\mu (x)<\infty, $$

by taking $\rho =\frac{1}{2M(C_{0}a)^{q}}$ in the last inequality, we get 
\begin{equation}
\frac{1}{2}\underset{0<s<N}{\sup }s^{q}\int_{E_{s}}w(x)d\mu (x)\leq \left(
1+\frac{1}{2(C_{0}a)^{q}}\right) \left( \frac{a}{\varepsilon }\right) ^{q}%
\underset{0<s<N\frac{\varepsilon }{a}}{\sup }s^{q}\int_{F_{s,1}}w(x)d\mu
(x).
\end{equation}
The desired inequality follows by letting $N$ goes to infinity.
\item[2)] Let $f$ be an arbitrary $\mu$-measurable function $f$. For any positive integer $k$, set $f_{k}=f\chi_{E_{k}}$ with \ $E_{k}=\left\{ x\in B_{\left( x_{0},k\right) }:\left| f(x)\right| \leq k\right\}.$  By part 1) of the proof, for any $k>0$, we have
\begin{equation}
\frac{1}{2}\underset{0<s<N}{\sup }s^{q}\int_{\left\{x\in X:I_{\gamma}f_{k}(x)>s\right\}}w(x)d\mu (x)\leq \left(
1+\frac{1}{2(C_{0}a)^{q}}\right) \left( \frac{a}{\varepsilon }\right) ^{q}%
\underset{0<s<N\frac{\varepsilon }{a}}{\sup }s^{q}\int_{\left\{x\in X:\mathcal M_{1,\frac{1}{\gamma}f_{k}(x)>s}\right\}}w(x)d\mu
(x).
\end{equation}
 So letting  $k$ goes to infinity, we obtain the result.
 \end{enumerate}
\epf

\begin{remark}
\label{TH: Continuité faib-Intfract} Assume that

$\bullet$ $\mu$ satisfies condition (\ref{normal}),

$\bullet$  $q$, $\theta $, $p$, $p_{1}$, $q_{1}$, $\theta _{1},$ $\gamma $  are elements of $ \left[ 0, \infty \right] $
such that  
\begin{equation}
1\leq q\leq \theta \leq p\ \text{with\ }0<\frac{1}{\theta }-\gamma =\frac{1
}{s}
\end{equation}
and 
\begin{equation}
q<q_{1}\leq \theta _{1}\leq p_{1}<\infty \ \text{with }0<\frac{1}{q_{1}}%
-\gamma =\frac{1}{t}\leq \frac{1}{p_{1}},
\end{equation}

$\bullet$ $\Phi $ is a doubling Young function whose conjugate function $\Phi ^{\ast
}$ satisfies the $B_{\frac{q_{1}}{q}}$ condition,

$\bullet$ $v$ and $w$ are
two weights for which there exists a constant $A$ such that  
\begin{equation*}
\mu \left( B\right) ^{-\frac{1}{t}}\left\| w\chi _{B}\right\|
_{t} \left\| v^{-q}\right\| _{\Phi ,B}^{\frac{1}{q}}
\leq A,\ \ \ B\text{ ball }
\end{equation*}
and $w^{t}$ satisfies $\mathcal{A}%
_{\infty }$ condition.

 Then, there is a constant $C$ such that for any
$\mu$-measurable functions $f$  and  $a>0$, we have : 
\begin{equation*}
\left( \int_{E_{a}}w^{t}(x)d\mu (t)\right) ^{\frac{1}{t}}\leq C\left(
a^{-1}\left\| fv\right\| _{q_{1},p_{1},\theta _{1}}\right) \left(
a^{-1}\left\| f\right\| _{q,p,\theta }\right) ^{s\left( \frac{1}{q_{1}}-%
\frac{1}{\theta _{1}}\right) },
\end{equation*}
where$\qquad E_{a}=\left\{ x\in X\text{ / }\left| I_{\gamma }f(x)\right|
>a\right\} $.
\end{remark}

\proof
Let $f$ be a $\mu$-measurable function. From Theorem \ref{TH: cont op max}, it follows that there exists a constant $C$ such that 
\begin{equation*}
\underset{a>0}{\sup }a^{1+s\left( \frac{1}{q_{1}}-\frac{1}{\theta _{1}}%
\right) }\left( \int_{E_{a}}w^{t}(x)d\mu (x)\right) ^{\frac{1}{t}}\leq C%
\underset{a>0}{\sup }a^{1+s\left( \frac{1}{q_{1}}-\frac{1}{\theta _{1}}%
\right) }\left( \int_{F_{a}}w^{t}(x)d\mu (x)\right) ^{\frac{1}{t}},
\end{equation*}
with $F_{a}=\left\{ x\in X:\mathcal M_{1,\frac{1}{\tau }}f(x)>a\right\} $.

Since $\mathcal M_{1,\beta }\leq \mathcal M_{q,\beta }$ for $q>1$, the result follows from Theorem \ref{TH: cont op max}.
\epf

\noindent{\bf Aknowledgment}We would like to tank S. Grellier for supplying us with documents. We are indebted to the refree for his valuable comments, and particularly for suggesting us not to impose normality condition on the underlying homogeneous space.

\end{document}